\documentclass[10pt,reqno]{amsart}
\usepackage{amsmath, hyperref}
\usepackage{amsthm}
\usepackage{amssymb}
\usepackage{amsfonts}
\usepackage{latexsym}
\usepackage{hyperref}
\usepackage{cite}
\usepackage{xcolor}

\usepackage{fancyvrb}
\theoremstyle{plain}
\newtheorem{theorem}{Theorem}

\newtheorem{alpthm}{Theorem}

\theoremstyle{definition}

\newtheorem{remark}[theorem]{Remark}

\newcommand{\comment}[1]{\marginpar{\color{blue} $\bullet$ \scriptsize #1}}
\renewcommand{\L}{\mathbb{L}}

\newcommand{\Q}{\mathbb{Q}}
\newcommand{\R}{\mathbb{R}}

\newcommand{\K}{\mathbb{K}}

\newcommand{\C}{\mathbb{C}}
\newcommand{\PP}{\mathfrak{p}}
\renewcommand{\O}{\mathcal{O}}

\newcommand{\Li}{\operatorname{Li}}
\newcommand{\Ei}{\operatorname{Ei}}

\renewcommand{\Re}{\operatorname{Re}}

\let\oldenumerate=\enumerate
	\def\enumerate{
	\oldenumerate
	\setlength{\itemsep}{5pt}
	}
\let\olditemize=\itemize
	\def\itemize{
	\olditemize
	\setlength{\itemsep}{5pt}
	}

\allowdisplaybreaks

\begin{document}
\VerbatimFootnotes
\title[Explicit Mertens' Theorems for Number Fields]{Explicit Mertens' Theorems for Number Fields}

	\author{Stephan Ramon Garcia}
	\address{Department of Mathematics, Pomona College, 610 N. College Ave., Claremont, CA 91711} 
	\email{stephan.garcia@pomona.edu}
	\urladdr{\url{http://pages.pomona.edu/~sg064747}}
	
\author{Ethan Simpson Lee}
\address{School of Science, UNSW Canberra at the Australian Defence Force Academy, Northcott Drive, Campbell, ACT 2612} 
\email{ethan.s.lee@student.adfa.edu.au}
\urladdr{\url{https://www.unsw.adfa.edu.au/our-people/mr-ethan-lee}}
	
\thanks{SRG supported by NSF Grant DMS-1800123.}



\begin{abstract}
Assuming the Generalized Riemann Hypothesis we obtain uniform, effective number-field analogues of Mertens' theorems.
\end{abstract}

\comment{To make the estimate \eqref{eq:MMK} below completely explicit requires a careful study of the Dedekind zeta residue $\kappa_{\K}$.
That long and difficult endeavor grew significantly and has since been spun off into a separate project \cite{DukePaper} (leaving the preprint you are currently reading in a temporarily truncated and diminished state).
We will eventually return to this preprint and complete it with improved and updated estimates, but \cite{DukePaper} and related projects that will contribute improved estimates for this project take priority. }

\maketitle
\section{Introduction}
Assuming the Generalized Riemann Hypothesis (GRH), we establish 
uniform, explicit number-field analogues of Mertens' classical theorems:
\begin{align}
\sum_{p\leq x} \frac{\log p}{p} &= \log x + O(1), \label{eq:1}\\
\sum_{p \leq x} \frac{1}{p} &= \log \log x + M + O\bigg(\frac{1}{\log x} \bigg), \label{eq:2}\\
\prod_{p\leq x} \bigg(1 -\frac{1}{p} \bigg) &= \frac{e^{-\gamma}}{\log x} \big(1+o(1) \big).\label{eq:3}
\end{align}
Here $p$ denotes a prime number, $M = 0.2614\ldots$ is the Meissel--Mertens constant, and $\gamma = 0.5772\ldots$ is the Euler--Mascheroni constant.  Mertens obtained these results in 1874 \cite{Mertens}, well before the prime number
theorem was proved independently by Hadamard \cite{hadamard1896distribution} and de la Vall\'ee Poussin \cite{de1896fonction}.
See Ingham \cite[Thm.~7]{Ingham} or Montgomery--Vaughan \cite[Thm.~2.7]{MontgomeryVaughan} for modern proofs of Mertens' theorems, and 
Rosen \cite[(3.17) - (3.30)]{Rosser} for explicit unconditional error bounds.

Rosen \cite[Lem.~2.3, Lem.~2.4, Thm.~2]{Rosen} generalized \eqref{eq:1}, \eqref{eq:2} and \eqref{eq:3} to the number-field setting without explicit error terms; see also \cite{Lebacque}. 
Assuming GRH, we provide explicit error bounds for these results.
Our results are uniform (valid for $x\geq 2$) and the constants involved depend only on the degree and discriminant of the number field,
and the corresponding Dedekind zeta residue.
This complements \cite{EMT4NF1}, in which the authors obtained similar results in the unconditional context.

\subsection*{Definitions}
Let $\K$ be a number field with ring of algebraic integers $\O_{\K}$. 
Let $n_{\K}$ denote the degree of $\K$ and $\Delta_{\K}$ the discriminant of $\K$. 
In what follows, $\PP\subset\O_{\K}$ denotes a prime ideal, $\mathfrak{a}\subset\O_{\K}$ an ideal, and $N(\mathfrak{a})$ the norm of $\mathfrak{a}$.

The Dedekind zeta function 
\begin{equation*}
\zeta_{\K}(s) = \sum_{\mathfrak{a} \subseteq \O_{\K}} \frac{1}{N(\mathfrak{a})^s}
= \prod_{\PP} \bigg(1 - \frac{1}{N(\PP)^s}\bigg)^{-1}
\end{equation*}
corresponding to $\K$ is analytic on $\C$ with exception of a simple pole at $s=1$. 
The analytic class number formula asserts that the residue of $\zeta_{\K}(s)$ at $s=1$ is
\begin{equation}\label{eqn:residue_class_ana_form}
\kappa_{\K} = \frac{2^{r_1}(2\pi)^{r_2}h_{\K}R_{\K}}{w_{\K}\sqrt{|\Delta_{\K}|}},
\end{equation}
where
$r_1$ is the number of real places of $\K$, 
$r_2$ is the number of complex places of $\K$, 
$w_{\K}$ is the number of roots of unity in $\K$, 
$h_{\K}$ is the class number of $\K$, 
and $R_{\K}$ is the regulator of $\K$ \cite{Lang}.
The nontrivial zeros of $\zeta_{\K}$ lie in the critical strip, $0 < \Re s < 1$, where there might exist an exceptional zero $\beta$, 
which is real and not too close to $\Re{s} =1$ \cite[p.~148]{Stark}. The Generalized Riemann Hypothesis (GRH) claims that 
the nontrivial zeros of $\zeta_{\K}(s)$ satisfy $\Re s = \frac{1}{2}$; in particular, $\beta$ does not exist.

Let $\K \subseteq \L$ be a Galois extension of number fields with Galois group $\mathcal{G}=\operatorname{Gal}(\L/\K)$. 
Let $n_{\L}$ denote the degree of $\L$ over $\K$ and let $\Delta_{\L}$ denote the absolute value of the discriminant of $\L$ over $\Q$.
Suppose that $\mathfrak{P}$ is a prime ideal of $\L$ lying above an unramified prime $\PP$
of $\K$.  The Artin symbol $[\L / \PP]$ denotes the conjugacy class of Frobenius automorphisms corresponding to prime ideals $\mathfrak{P}|\PP$.
For each conjugacy class $\mathcal{C}\subset \mathcal{G}$, the prime ideal counting function is
\begin{equation*}
\pi_{\mathcal{C}}(x, \L/\K)\, = \,\#\left\{\PP : \text{$\PP$ is unramified in $\L$},\,\, \left[\tfrac{\L/\K}{\PP}\right] = \mathcal{C},\,\, N_{\K}(\PP)\leq x\right\},
\end{equation*}
in which $N_{\K}(\cdot)$ denotes the norm in $\K$.
If $\K = \L$, then we may specialize our notation so that the prime ideal counting function is 
$\pi_{\K}(x) = \sum_{N(\PP)\leq x} 1$.

\subsection*{Statement of results}
Assuming GRH, we obtain the following number-field analogues of
Mertens' theorems \eqref{eq:1}, \eqref{eq:2} and \eqref{eq:3}.
The proof involves a recent estimate of Greni\'e--Molteni \cite{Grenie2019} for $\pi_{\K}(x)$.

\begin{alpthm}\label{Theorem:Main}
Assume GRH.
For a number field $\K$ and $x\geq 2$,
\begin{align}
    \sum_{N(\PP)\leq x}\frac{\log{N(\PP)}}{N(\PP)} &\,=\, \log{x} + A_{\K}(x) ,\label{eq:A1}\tag{A1}\\[5pt]
    \sum_{N(\PP)\leq x}\frac{1}{N(\PP)} &\,=\, \log\log{x} + M_{\K} + B_{\K}(x),\label{eq:B1} \tag{B1}\\[5pt]
    \prod_{N(\PP)\leq x}\left(1 - \frac{1}{N(\PP)}\right) &\,=\, \frac{e^{-\gamma}}{\kappa_{\K}\log{x}} \big(1 + C_{\K}(x) \big),\label{eq:C1}\tag{C1}
\end{align}
in which
\begin{align}
M_{\K} &=\, \gamma + \log{\kappa_{\K}} + \sum_{\PP}\left[\frac{1}{N(\PP)} + \log\left(1 - \frac{1}{N(\PP)}\right)\right],\nonumber\\
|A_{\K}(x)|  &\,\leq\,  4.73 \log |\Delta_{\K}| + 9.27 \, n_{\K} + \log 2,  \label{eq:A2}\tag{A2}\\[5pt]
|B_{\K}(x)| &\,\leq\,  \frac{13.47 \log |\Delta_{\K}|+ (26.37 + 0.12 \log x) \, n_{\K}}{\sqrt{x}}, \qquad \text{and}
\label{eq:B2}\tag{B2}\\[5pt]
|C_{\K}(x)| &\,\leq\, |E_{\K}(x)|e^{| E_{\K}(x)|} \qquad \text{with}\qquad |E_{\K}(x)| \,\leq\, \frac{n_{\K}}{x-1} + |B_{\K}(x)|.\label{eq:C2}\tag{C2}
\end{align}
In particular, $E_{\K}(x) = o(1)$ (hence $C_{\K}(x) = o(1)$ as $x \to \infty$) and
\begin{equation}\label{eq:MMK}\tag{M}
\gamma + \log \kappa_{\K} -n_{\K} \,\,\leq\,\, M_{\K} \,\,\leq\,\,  \gamma + \log \kappa_{\K}.
\end{equation}
\comment{Earlier incarnations of this preprint included explicit estimates of the Dedekind zeta residue $\kappa_{\K}$
which have since been migrated to \cite{DukePaper}.  We plan to return to the current project when our investigations into Dedekind zeta residue 
and related topics are completed.}%
\end{alpthm}
\vspace{-10pt}

For $n_{\K} = 1$ ($\K = \Q$), Schoenfeld \cite[Cor.~2-3]{Schoenfeld} provides smaller constants.
In particular, $\kappa_{\Q} = 1$ so we can restrict our attention to $n_{\K} \geq 2$.
Although unconditional bounds for $\kappa_{\K}$ are available (see
\cite{Louboutin00, Louboutin01, Louboutin03, Louboutin11, Stark} and the discussion in \cite{EMT4NF1}),
stronger bounds are possible under GRH and the assumption that $\zeta_{\K}(s) / \zeta(s)$ is entire \cite{DukePaper}.

Slight changes to the proof of Theorem \ref{Theorem:Main} (see Remark \ref{Remark:Slight}) yield the following.

\begin{alpthm}\label{Theorem:MainExtended}
Assume GRH and let $\K \subseteq \L$ be a Galois extension of number fields. For $x \geq 2$,
\begin{align*}
    \sum_{N(\PP)\leq x}\frac{\log{N(\PP)}}{N(\PP)} &\,=\, \frac{\# \mathcal{C}}{\# \mathcal{G}}\log{x} + A_{\L}(x) ,\\
    \sum_{N(\PP)\leq x}\frac{1}{N(\PP)} &\,=\, \frac{\# \mathcal{C}}{\# \mathcal{G}} \log\log{x} + M_{\L} + B_{\L}(x),\\
    \prod_{N(\PP)\leq x}\left(1 - \frac{1}{N(\PP)}\right) &\,=\, \frac{e^{-\gamma}}{\kappa_{\L}(\log{x})^{\# \mathcal{C}/\# \mathcal{G}}} \big(1 + C_{\L}(x) \big),
\end{align*}
in which $\PP$ is a prime ideal in $\K$ that does not ramify in $\L$,
$\kappa_{\L}$ denotes the residue of $\zeta_{\L}(s)$ at $s=1$, and $A_{\L}(x)$, $B_{\L}(x)$, $C_{\L}(x)$, $M_{\L}$ are defined as $A_{\K}(x)$, $B_{\K}(x)$, $C_{\K}(x)$, $M_{\K}$ are in Theorem \ref{Theorem:Main}, but with $n_{\L}$, $\Delta_{\L}$, $\kappa_{\L}$ in place of $n_{\K}$, $|\Delta_{\K}|$, $\kappa_{\K}$.
Our explicit bounds for $M_{\K}$ and $\kappa_{\K}$ also carry forward for $M_{\L}$ and $\kappa_{\L}$.
\end{alpthm}


\subsection*{Acknowledgments}%
\comment{Most of these acknowledgments are for comments on material that has migrated to \cite{DukePaper}.}%
We thank Karim Belabas, Eduardo Friedman, Lenny Fukshansky, Peter Kim, St\'ephane Louboutin, and Tim Trudgian for helpful feedback.

\section{Proof of Theorem \ref{Theorem:Main}}\label{Section:Proof}
We prove \eqref{eq:A1}, \eqref{eq:B1}, and \eqref{eq:C1}, along with the associated error bounds \eqref{eq:A2}, \eqref{eq:B2}, and \eqref{eq:C2}, respectively, in separate subsections. We then prove \eqref{eq:MMK}, which is a consequence of the computations in the proof of \eqref{eq:C1}.  

\subsection{Preliminaries}
We require an explicit version of the prime ideal theorem, which was originally established by Landau in 1903 
(without explicit constants) \cite{Landau}. 
Recall that $ \pi_{\K}(x) = \sum_{N(\PP)\leq x} 1$ is the prime-ideal counting function and let
\begin{equation*}
    \Li(x) = \int_2^x \frac{dt}{\log{t}}
\end{equation*}
denote the offset logarithmic integral. 
Under GRH, 
\begin{equation}\label{eq:PITstatement}
    \pi_{\K}(x) = \Li(x) + R_{\K}(x),
\end{equation}
in which $R_{\K}(x) = O(x^{1/2}\log{x})$.
To bound $R_{\K}(x)$ explicitly, one can follow Lagarias--Odlyzko \cite{LagariasOdlyzko}.
Since the quality of the error term one obtains depends upon the zero-free region of $\zeta_{\K}$ one uses,
GRH naturally enters the conversation.
In 2019, Greni\'e--Molteni \cite[Cor.~1]{Grenie2019} proved that GRH implies that \eqref{eq:PITstatement}
holds for $x \geq 2$ with
\begin{equation}\label{eq:Grenie}
|R_{\K}(x)|
\leq \sqrt{x}\left[\left(\frac{1}{2\pi} + \frac{3}{\log x}\right)\log|\Delta_{\K}|
 + \left(\frac{\log x}{8\pi} + \frac{1}{4\pi} + \frac{6}{\log x} \right) n_{\K}\right].
\end{equation}
This strengthens a result announced by Oesterl\'{e} \cite{Oesterle} and improves the constants Winckler obtains \cite[Thm.~1.2]{Winckler} by using methods from Schoenfeld \cite{Schoenfeld}, Lagarias--Odlyzko \cite{LagariasOdlyzko}, and previous work of Greni\'e--Molteni \cite{Grenie2016}.  

The final important ingredient we need is the following formula for a continuously-differentiable function $f$. 
Integration by parts and \eqref{eq:PITstatement} provide
\begin{align}
    \sum_{N(\PP)\leq x} f(N(\PP))
    &= f(x) \pi_{\K}(x) - \int_2^x f'(t)\pi_{\K}(t)\,dt \nonumber \\
    &=f(x) \Li(x)  - \int_2^x f'(t) \Li(t)\,dt+ f(x) R_{\K}(x)  - \int_2^x f'(t) R_{K}(t)\,dt\nonumber \\
    &= \int_2^x \frac{f(t)}{\log{t}}dt + f(x)R_{\K}(x) - \int_2^x f'(t)R_{\K}(t)dt. \label{eq:IBP2}
\end{align}

\begin{remark}\label{Remark:Slight}
We briefly remark how to adjust the proof below to obtain Theorem \ref{Theorem:MainExtended}.
For $x \geq 2$, Greni\'e--Molteni \cite[Cor.~1]{Grenie2019} proved 
\begin{equation*}
\pi_{\mathcal{C}}(x) =  \frac{\# \mathcal{C}}{\# \mathcal{G}} \Li(x) + R_{\L/\K}(x),
\end{equation*}
in which 
\begin{equation*}
|R_{\L/\K}(x)|
\leq\frac{\# \mathcal{C}}{\# \mathcal{G}} \sqrt{x}\left[\left(\frac{1}{2\pi} + \frac{3}{\log x}\right)\log\Delta_{\L} + \left(\frac{\log x}{8\pi} + \frac{1}{4\pi} + \frac{6}{\log x} \right) n_{\L}\right].
\end{equation*}
Since $\# \mathcal{C} / \# \mathcal{G} \leq 1$, we can bound $|R_{\L/\K}(x)|$ above
by \eqref{eq:Grenie}.  The proof below goes through \emph{mutatis mutandis}.
For example, any occurrence of $\Li(x)$
should be replaced by  $(\# \mathcal{C}/\# \mathcal{G}) \Li(x)$.
\end{remark}

\subsection{Proof of \eqref{eq:A1} and \eqref{eq:A2}}
Let $f(t) = \log{t}/t$ in \eqref{eq:IBP2} and obtain
\begin{equation}\label{eq:AKGRH}
\sum_{N(\PP)\leq x} \frac{\log{N(\PP)}}{N(\PP)}
    = \log{x} + \underbrace{f(x)R_{\K}(x) - \int_2^x f'(t)R_{\K}(t)dt - \log{2}}_{A_{\K}(x)}.
\end{equation}
For $x \geq 2$, $\log x / \sqrt{x}$ achieves its maximum value $2/e$ at $x= e^2$ and $(\log x)^2 / \sqrt{x}$ achieves its maximum value $16/e^2$ at $x= e^4$.
Thus, \eqref{eq:Grenie} yields
\begin{align}
    |f(x)R_{\K}(x)|
    &\leq \frac{\log{x}}{\sqrt{x}}\left[\left(\frac{1}{2\pi} + \frac{3}{\log x}\right)\log|\Delta_{\K}| + \left(\frac{\log x}{8\pi} + \frac{1}{4\pi} + \frac{6}{\log x} \right) n_{\K}\right]. \nonumber\\
    &= \frac{1}{\sqrt{x}}\left(3\log|\Delta_{\K}| + 6 n_{\K}\right)
    + \frac{\log{x}}{\sqrt{x}}\left(\frac{\log|\Delta_{\K}|}{2\pi} + \frac{n_{\K}}{4\pi} \right)
    + \frac{(\log x)^2}{\sqrt{x}}\frac{n_{\K}}{8\pi} \nonumber \\
    &\leq \frac{1}{\sqrt{2}}\left(3\log|\Delta_{\K}| + 6 n_{\K}\right)
    + \frac{2}{e}\left(\frac{\log|\Delta_{\K}|}{2\pi} + \frac{n_{\K}}{4\pi} \right)
    + \frac{2n_{\K}}{e^2\pi} \nonumber \\
    &<2.23843 \log |\Delta_{\K}| + 4.3874\, n_{\K}. \label{eq:AddThis1}
\end{align}

Observe that
\begin{equation*}
f'(t)= \frac{1 - \log t}{t^2} 
\end{equation*}
undergoes a sign change from positive to negative at $x=e$ and that
\begin{equation*}
\frac{\log x}{8\pi} + \frac{1}{4\pi} + \frac{6}{\log x}, \quad \text{whose derivative is}\quad  \frac{1}{8 \pi  x}-\frac{6}{x (\log x)^2},
\end{equation*}
decreases until $e^{4\sqrt{3 \pi}} \approx 215{,}328.6$ (in particular, it is decreasing on $[2,e]$).
Consequently, we split the integral in \eqref{eq:AKGRH} at $x = e$. First note that
\begin{align}
&\left| \int_2^e f'(t) R_{\K}(t)\,dt \right| \nonumber \\
&\quad\leq \int_2^e \frac{1 - \log t}{t^{3/2}}
 \left[\left(\frac{1}{2\pi} + \frac{3}{\log t}\right)\log|\Delta_{\K}|
 + \left(\frac{\log t}{8\pi} + \frac{1}{4\pi} + \frac{6}{\log t} \right) n_{\K}\right] dt \nonumber \\
&\quad \leq \left(\frac{4}{\sqrt{e}}-\sqrt{2} (1+\log 2)\right)\left[\left(\frac{1}{2\pi} + \frac{3}{\log 2}\right)\log|\Delta_{\K}|
 + \left(\frac{\log 2}{8\pi} + \frac{1}{4\pi} + \frac{6}{\log 2} \right) n_{\K}\right] \nonumber \\
&\quad< 0.14203 \log |\Delta_{\K}| + 0.27737 \,n_{\K}. \label{eq:AddThis2}
\end{align}
For $x \geq e$, symbolic integration provides
\begin{align}
&\left|  \int_e^x f'(t) R_{\K}(t)\,dt \right| \nonumber \\
&\qquad\leq \int_e^{\infty}  \frac{\log t - 1 }{t^{3/2}}\left[\left(\frac{1}{2\pi} + \frac{3}{\log t}\right)\log|\Delta_{\K}|
 + \left(\frac{\log t}{8\pi} + \frac{1}{4\pi} + \frac{6}{\log t} \right) n_{\K}\right] dt \nonumber \\
&\qquad= \frac{(3 \pi  \sqrt{e} \Ei (-\frac{1}{2})+6 \pi +2) \log |\Delta_{\K}| }{\pi\sqrt{e}  }+\frac{(12 \pi  \sqrt{e} \Ei(-\frac{1}{2})+24 \pi +7) n_{\K}}{2\pi \sqrt{e}  } \nonumber \\
&\qquad < 2.34600 \log |\Delta_{\K}| + 4.59546 \,n_{\K} , \label{eq:AddThis3}
\end{align}
in which $\Ei(x)$ is the exponential integral.  

In light of \eqref{eq:AKGRH}, we add \eqref{eq:AddThis1}, \eqref{eq:AddThis2}, and \eqref{eq:AddThis3}, and obtain the bound
\begin{equation*}
|A_{\K}(x)| <   4.72646 \log |\Delta_{\K}| + 9.26023 \, n_{\K} + \log 2.
\end{equation*}
This concludes the proofs of \eqref{eq:A1} and \eqref{eq:A2}. \qed

\subsection{Proof of \eqref{eq:B1} and \eqref{eq:B2}}

Let $f(t) = 1/t$ in \eqref{eq:IBP2} and obtain
\begin{align*}
\sum_{N(\PP)\leq x} \frac{1}{N(\PP)}
&=\int_2^x \frac{dt}{t\log{t}} + \frac{R_{\K}(x)}{x} + \int_2^x\frac{R_{\K}(t)\,dt}{t^2} \\
&= \log \log x - \log \log 2 + \frac{R_{\K}(x)}{x} + \int_2^{\infty}\frac{R_{\K}(t)\,dt}{t^2} - \int_x^{\infty}\frac{R_{\K}(t)\,dt}{t^2} \\
&= \log \log x   + \underbrace{ \int_2^{\infty}\frac{R_{\K}(t)}{t^2}dt - \log \log 2 }_{M_{\K}}
+ \underbrace{ \frac{R_{\K}(x)}{x}- \int_x^{\infty}\frac{R_{\K}(t)}{t^2}dt}_{B_{\K}(x)},
\end{align*}
in which the convergence of the improper integral is guaranteed by \eqref{eq:Grenie}.

For $x\geq 2$, \eqref{eq:Grenie} provides
\begin{equation*}
\frac{|R_{\K}(x)|}{x} 
\leq \frac{1}{\sqrt{x}}\left[\left(\frac{1}{2\pi} + \frac{3}{\log x}\right)\log|\Delta_{\K}|
 + \left(\frac{\log x}{8\pi} + \frac{1}{4\pi} + \frac{6}{\log x} \right) n_{\K}\right] 
\end{equation*}
and
\begin{align*}
&\int_x^{\infty}\frac{| R_{\K}(t)| \,dt}{t^2} \\
    &\qquad\leq \int_x^\infty \frac{1}{t^{3/2}}\left[\left(\frac{1}{2\pi} + \frac{3}{\log t}\right)\log|\Delta_{\K}| + \left(\frac{\log t}{8\pi} + \frac{1}{4\pi} + \frac{6}{\log t} \right) n_{\K}\right]dt \\
    &\qquad= \log|\Delta_{\K}| \int_x^\infty \frac{1}{t^{3/2}}\left(\frac{1}{2\pi} + \frac{3}{\log t}\right)dt + n_{\K}\int_x^\infty \frac{1}{t^{3/2}} \left(\frac{\log t}{8\pi} + \frac{1}{4\pi} + \frac{6}{\log t} \right) dt\\
    &\qquad\leq \log|\Delta_{\K}| \int_x^\infty \frac{1}{t^{3/2}}\left(\frac{1}{2\pi} + \frac{3}{\log x}\right)dt + n_{\K}\int_x^\infty \frac{1}{t^{3/2}} \left(\frac{\log t}{8\pi} + \frac{1}{4\pi} + \frac{6}{\log x} \right) dt\\
&\qquad=  \frac{1}{\sqrt{x}} \left[\left( \frac{1}{\pi}+\frac{6}{ \log x} \right) \log |\Delta_{\K}| + 
\left(\frac{\log x}{4 \pi }+\frac{1}{\pi }+\frac{12}{\log x}\right) n_{\K} \right].
\end{align*}
For $x \geq 2$, the triangle inequality provides 
\begin{align}
|B_{\K}(x)|
&\leq \frac{1}{\sqrt{x}} \left[ \left(\frac{3}{2 \pi } + \frac{9}{\log 2}\right) \log |\Delta_{\K}| + \left(\frac{3 \log x}{8 \pi }+\frac{5}{4 \pi }+\frac{18}{\log 2}\right) n_{\K} \right] \label{eq:UseThisLater} \\
&< \frac{1}{\sqrt{x}}\left[ 13.47 \log |\Delta_{\K}|+ (26.37 + 0.120 \log x)  \,n_{\K} \right].\nonumber
\end{align}

Now we find the constant $M_{\K}$, following Ingham \cite{Ingham}. Since this part of the proof is identical to that in \cite{EMT4NF1}, we present an abbreviated version here.

For $\Re s > 1$, partial summation and \eqref{eq:B1} provide
\begin{align}
\sum_{\PP}&\frac{1}{N(\PP)^s}
= \lim_{x\to\infty} \bigg( \sum_{N(\PP)\leq x}\frac{1}{N(\PP)^s} \bigg)
= \lim_{x\to\infty} \bigg( \sum_{N(\PP)\leq x}\frac{1}{N(\PP)^{s-1} N(\PP)} \bigg) \nonumber\\
&= \lim_{x\to\infty} \bigg(\frac{1}{x^{s - 1}}\sum_{N(\PP)\leq x}\frac{1}{N(\PP)}\bigg) + (s-1) \int_{2}^{\infty}\bigg(\sum_{N(\PP)\leq t}\frac{1}{N(\PP)}\bigg) \frac{dt}{t^{s}}\nonumber\\
&=(s-1) \int_{2}^{\infty}\bigg(\sum_{N(\PP)\leq t}\frac{1}{N(\PP)}\bigg) \frac{dt}{t^{s}}\nonumber\\
&=   (s-1)\int_{2}^\infty\frac{M_{\K}}{t^{s}}\,dt
+ (s-1)\int_{2}^\infty\frac{B_{\K}(t)}{t^{s}} \,dt 
+ (s-1)\int_{2}^\infty\frac{\log\log t}{t^{s}}\,dt.\label{eqn:refmel8er}
\end{align}
As $s \to 1^+$,  the first integral tends to $M_{\K}$, the second tends to zero by \eqref{eq:B2},
and the substitution $t^{s-1} = e^y$ converts the third integral into
\begin{equation*}
(s-1)\int_{2}^{\infty}\frac{\log\log{t}}{t^{s}}\,dt
= \int_{\log(2^{s-1})}^{\infty} e^{-y} \log y \, dy - 2^{1-s} \log(s-1),
\end{equation*}
which tends to $- \gamma- \log(s-1)$.  As $s \to 1^+$, the Euler product formula yields
\begin{align}
M_{\K} &-\gamma +o(1)\nonumber\\
&=\log(s-1) + \sum_{\PP}\left[ \frac{1}{N(\PP)^s} + \log\left(1 - \frac{1}{N(\PP)^s}\right) \right] - \sum_{\PP}  \log\left(1 - \frac{1}{N(\PP)^s}\right) \nonumber \\
&= \log\big((s-1)\zeta_{\K}(s)\big) + \sum_{\PP}\left[ \frac{1}{N(\PP)^s} + \log\left(1 - \frac{1}{N(\PP)^s}\right) \right] \nonumber \\
&= \log \kappa_{\K} + \sum_{\PP}\left[\frac{1}{N(\PP)} + \log\left(1 - \frac{1}{N(\PP)}\right)\right] +o(1), \label{eq:MKGo1}
\end{align}
in which the last sum converges, by comparison with $\sum_{\PP} N(\PP)^{-2}$.   \qed


\subsection{Proofs of \eqref{eq:C1} and \eqref{eq:C2}}

The proof is identical to that in \cite{EMT4NF1} so we just sketch it here for the sake of completeness.
From \eqref{eq:MKGo1}, we may write
\begin{equation}\label{eq:refmelastX}
    -\gamma - \log{\kappa_{\K}} + M_{\K} 
    = \sum_{N(\PP) \leq x}\left[\frac{1}{N(\PP)} + \log\left(1 - \frac{1}{N(\PP)}\right)\right] + F_{\K}(x).
\end{equation}
Since 
\begin{equation}\label{eq:yyy}
0 \leq -y -\log(1-y) \leq  \frac{y^2}{1 - y},
\end{equation}
it follows that
\begin{align}
    | F_{\K}(x)|
    &=- \sum_{N(\PP) > x}\left[\frac{1}{N(\PP)} + \log\left(1 - \frac{1}{N(\PP)}\right)\right] \nonumber \\
    &\leq \sum_{N(\PP) > x}\frac{1}{N(\PP) (N(\PP) - 1)} \nonumber \\
    &\leq \sum_{p > x}\sum_{f_i}\frac{1}{p^{f_i} (p^{f_i} - 1)} 
    < n_{\K}\sum_{m > x}\frac{1}{m (m - 1)}  
    \leq \frac{n_{\K}}{ x-1 }, \label{eq:FKx}
\end{align}
in which $\sum_{f_i}$ denotes the sum over the inertia degrees $f_i$ of the prime ideals lying
over $\PP$. 
Then \eqref{eq:B1} and \eqref{eq:refmelastX} imply
\begin{equation*}
-\gamma - \log{\kappa_{\K}} + M_{\K} =  \sum_{N(\PP) \leq x} \log\left(1 - \frac{1}{N(\PP)}\right) + \log\log{x} + M_{\K} +
\underbrace{B_{\K}(x) + F_{\K}(x)}_{E_{\K}(x)}.
\end{equation*}
Exponentiate and simplify to get
\begin{equation*}
    \prod_{N(\PP) \leq x} \left(1 - \frac{1}{N(\PP)}\right) 
    = \frac{e^{- \gamma}}{\kappa_{\K} \log{x}} e^{-E_{\K}(x)}
    = \frac{e^{- \gamma}}{\kappa_{\K} \log{x}}\big(1 + C_{\K}(x) \big),
\end{equation*}
where $|C_{\K}(x)| \leq  |E_{\K}(x)| e^{ |E_{\K}(x)|}$ since $|e^t - 1|  \leq |t| e^{|t|}$ for $t \in \R$. \qed

\subsection{Proof of \eqref{eq:MMK}}
Since there are no prime ideals of norm less than $2$,
\eqref{eq:MKGo1} yields
\begin{equation*}
M_{\K} = \gamma + \log\kappa_{\K} + F_{\K}(2-\delta)\quad \text{for $\delta \in (0,1)$},
\end{equation*}
in which $F_{\K}(x)$ is defined in \eqref{eq:refmelastX} 
In particular, \eqref{eq:FKx} reveals that
\begin{equation*}
    |F_{\K}(2-\delta)| \leq n_{\K} \qquad\text{as}\qquad\delta\to 0^+.
\end{equation*}
In light of \eqref{eq:yyy}, each summand in $F_{\K}(2-\delta)$ is non-positive,
so $F_{\K}(2-\delta) \leq 0$ as $\delta\to0^+$. Consequently, $- n_{\K} \leq F_{\K}(2-\delta) \leq 0$ and hence
\begin{equation*}
-n_{\K} \leq M_{\K} - \gamma - \log \kappa_{\K} \leq 0,
\end{equation*}
which is equivalent to \eqref{eq:MMK}. \qed

\bibliographystyle{amsplain}
\bibliography{EMTNF-GRH}

\end{document}